\documentclass[12pt]{article}
\usepackage{amssymb,amsthm,amsmath}
\usepackage{enumerate}
\usepackage{verbatim}
\usepackage{cite}

\usepackage{graphicx}

\newcommand{\Var}{\mathop{\rm Var}}

\def\csect#1{Section~\ref{#1}}

\let\Prp=\Pr \def\Pr{\Prp\nolimits}
\newcommand{\be}{\begin{equation}}
\newcommand{\ee}{\end{equation}}

\newtheorem{theorem}{Theorem}
\newtheorem{corollary}[theorem]{Corollary}

{\theoremstyle{definition}
\newtheorem*{conjecture}{Fluctuation Conjecture}
 
\newtheorem{remark}{Remark}
\newtheorem{definition}{Definition} }

\newcommand{\cthm}[1]{Theorem~\ref{#1}}

\newcommand{\cdef}[1]{Definition~\ref{#1}}

\newcommand{\ccor}[1]{Corollary~\ref{#1}}

\newcommand{\crem}[1]{Remark~\ref{#1}}
\newcommand{\cconjn}{the Fluctuation Conjecture}
\newcommand{\cconj}{the Fluctuation Conjecture }
\newcommand{\Cconj}{The Fluctuation Conjecture }

    \def\bbz{\mathbb{Z}}

\def\Cov{\mathop{\rm Cov}}

\def\sss{\scriptscriptstyle}
\def\interiorM{\rlap{\raise9pt\hbox to9.5pt{\hss$\sss\circ$}}M}

\def\etal{\eta^{(\Lt)}}\def\etaldh{\hat\eta^{(\Lt,\delta)}}
\def\etalzh{\hat\eta^{(\Lt,0)}}

\long\def\kill#1\endkill{\relax}

\def\0{{\it0}}\def\1{{\it1}}\def\2{{\it2}}\def\3{{\it3}}

\def\Nd{N^{(\delta)}}\def\Ndr{N^{(\delta)}_{\rm ren}}

\def\Ydh{\widehat Y^{(\delta)}}\def\Yzh{\widehat Y^{(0)}}
\def\Ndrh{\widehat N^{(\delta)}_{\rm ren}}
\def\Ndlh#1{\widehat N^{(\delta)}_{\le #1}}\def\Nzlh#1{\widehat N^{(0)}_{\le #1}}
\def\Ndgh#1{\widehat N^{(\delta)}_{> #1}}\def\Nzgh#1{\widehat N^{(0)}_{> #1}}
\def\Nzr{N^{(0)}_{\rm ren}}\def\Nzrh{\widehat N^{(0)}_{\rm ren}}

\def\Xdh{\widehat X^{(\delta)}}\def\Xzh{\widehat X^{(0)}}
\def\Fd{F^{(\delta)}}\def\sd{\sigma^{(\delta)}}
\def\pd{p_\delta}
\def\qd{q_\delta}\def\qz{q_0}\def\rd{r_\delta}

\def\IL{J_L}\def\ILt{J_{\Lt}}
\def\gtr{g^T_\rho}

\def\bll{\noindent\hbox to 15pt{\hfil$\bullet$\hfil}} 

\def\Lt{\widetilde L}

\newcommand\latop[2]{\genfrac{}{}{0pt}{}{#1}{#2}}

\def\murho{\mu^{(\rho)}}\def\nurho{\nu_\rho}\def\nurhoc{\nu_{\rho_c}}

\begin{document} \title{Approach to Hyperuniformity in the
  One-Dimensional Facilitated Exclusion Process}

\author{S. Goldstein\footnote{Department of Mathematics,
Rutgers University, New Brunswick, NJ 08903.},
J. L. Lebowitz\footnotemark[1],$\!$
\footnote{Also Department of Physics, Rutgers.}\ \ 
and E. R. Speer\footnotemark[1]\ \footnote{Corresponding author: speer@math.rtugers.edu}}
\date{June 11, 2025}

\maketitle

\begin{flushleft}\noindent {\bf Keywords:} Facilitated exclusion
  process, conserved lattice gas, hyperuniform states, approach to
hyperuniformity.

\par\medskip\noindent
 {\bf AMS subject classifications:} 60K35, 60K10, 82C22, 82C23

\end{flushleft}

\begin{abstract} For the one-dimensional Facilitated Exclusion Process
with initial state a product measure of density $\rho=1/2-\delta$,
$\delta\ge0$, there exists an infinite-time limiting state $\nu_\rho$
in which all particles are isolated and hence cannot move.  We study
the variance $V(L)$, under $\nu_\rho$, of the number of particles in
an interval of $L$ sites.  Under $\nu_{1/2}$ either all odd or all
even sites are occupied, so that $V(L)=0$ for $L$ even and $V(L)=1/4$
for $L$ odd: the state is {\it hyperuniform} \cite{T}, since $V(L)$ grows more
slowly than $L$.  We prove that for densities approaching 1/2 from
below there exist three regimes in $L$, in which the variance grows at
different rates: for $L\gg\delta^{-2}$, $V(L)\simeq\rho(1-\rho)L$,
just as in the initial state; for $A(\delta)\ll L\ll\delta^{-2}$, with
$A(\delta)=\delta^{-2/3}$ for $L$ odd and $A(\delta)=1$ for $L$ even,
$V(L)\simeq CL^{3/2}$ with $C=2\sqrt{2/\pi}/3$; and for
$L\ll\delta^{-2/3}$ with $L$ odd, $V(L)\simeq1/4$.  The analysis is
based on a careful study of a renewal process with a long tail.  Our
study is motivated by simulation results showing similar behavior in
higher dimensions; we discuss this background briefly.  \end{abstract}

\section{Introduction\label{intro}}

In the Facilitated Exclusion Process (FEP) on $\bbz^d$, also known as
the Conserved Lattice Gas process, each site of the lattice can be
occupied by at most one particle, so that a configuration $\eta$ is an
element of the configuration space $X=\{0,1\}^{\bbz^d}$.  At
Poisson-distributed times a particle chooses one of its
nearest-neighbor sites at random and attempts to jump to it,
succeeding only if the target site is unoccupied and the original site
has at least one occupied ({\it facilitating}) neighbor.  (Variations
of this dynamics, with simultaneous updating or with some other rule
for choosing the target site, have also been considered \cite{GLS4}.)
Note that when $d=1$ a chosen particle can jump in at most one
direction.  We will always assume that the system is started in a
Bernoulli initial state $\murho_0$ of density $\rho<1$, that is, a
product measure in which the $\eta_i$ are independent and take value
1 with probability $\rho$.  The evolved state at time $t$ will then be
denoted $\murho_t$; it is clearly translation invariant (TI).

The evolution of this system, or of minor variations of it, has been
investigated for $d=1$
\cite{AGLS,BBCS,BESS,BES,GKR,GR,GLS1,ZC,GLS2,GLS3,GLS4,O}, primarily
theoretically, and for $d\ge2$ \cite{GLS4,HL,L,RPV}, primarily via
simulation in a cubical box with periodic boundary conditions.  These
investigations strongly suggest the existence of a TI limiting state
 \be\label{nurho}
\nurho:=\lim_{t\to\infty}\murho_t.
 \ee
 This existence was proved for $d=1$ and $\rho\le1/2$ in \cite{ZC}
(and the limiting state was first described there).  Moreover, there
appears to be a {\it critical density} $\rho_c$ such that if
$\rho\le\rho_c$ then $\nurho$ is a {\it frozen} state in which all
particles are isolated and hence unable to move, while if
$\rho>\rho_c$ then $\nurho$ is an {\it active} stationary state in
which there is a nonzero density of particles with an occupied
neighboring site.  Necessarily $\rho_c\le1/2$, since for $\rho>1/2$ it
is geometrically impossible for all the particles to be isolated, and
indeed equality holds for $d=1$.  But for $d\ge2$, simulations suggest
values of $\rho_c$ which are much smaller, for example,
$\rho_c\approx0.3308$ \cite{GLS4}  for $d=2$.

Our main interest here will be in {\it fluctuations} in the measure
$\nurho$, that is, in the variance $V_\rho(L):=\Var_{\nurho}(N(L))$ of
the number $N(L)$ of particles in a cubical box of side $L$.  In
general, if such a variance computed from a TI measure $\mu$ grows as
$L^d$ when $L\nearrow\infty$, we say that $\mu$ has {\it normal
fluctuations}.  Hexner and Levine \cite{HL} observe that, in 2d and
3d, $\nurhoc$ is not of this character but rather is {\it
hyperuniform} \cite{T}: $V_{\rho_c}(L)$ grows more slowly than $L^d$,
specifically, $V_{\rho_c}(L)\sim L^{\lambda_1}$, with
$\lambda_1\approx1.57$ in 2d and $\lambda_1\approx2.76$ in 3d
\cite{HL}.  (We will relate two quantities by ``$\sim$'' or,
respectively, by ``$\simeq$'', to express the fact that asymptotically
their ratio is bounded away from both 0 and $\infty$ or, respectively,
is equal to 1.  The symbol ``$\simeq$'', used extensively in this
paper to express and prove our results, will be more precisely defined
in \cdef{asymp}.)

Hexner and Levine also discuss the behavior of $V_\rho(L)$ as
$\rho\nearrow\rho_c$.  Further consideration of this behavior has led
us \cite{GLS4} to the following conjecture, which we state for
$d\ge2$; the $d=1$ version is \cthm{main} below. We introduce the
notation $\delta=\rho_c-\rho>0$.

\begin{conjecture} For the FEP with $d\ge2$ a critical density
$\rho_c$ as described above exists, and $\nurhoc$ is hyperuniform.
For $\rho$ less than but close to $\rho_c$, three regimes in $L$ may
be identified.  For {\it small} $L$ (but still with $L\gg1$) the
variances grow approximately as in the hyperuniform state at $\rho_c$:
$V_\rho(L)\simeq C_1L^{\lambda_1}$. At some (approximately defined)
scale $L_1(\delta)$ the variances enter the regime of {\it
intermediate} $L$, in which they grow as
$V_\rho(L)\simeq C_2(\delta)L^{\lambda_2}$ with
$\lambda_2>d>\lambda_1$ and $C_2(\delta)>0$. Then above an
(approximate) scale $L_2=L_2(\delta)$ the growth is as
$V_\rho(L)\simeq\rho(1-\rho)L^d$, that is, exactly as in the initial
Bernoulli measure $\murho_0$. Finally, as $\rho\nearrow\rho_c$,
$L_1(\delta)$ and $L_2(\delta)$ increase as
$L_i\sim \delta^{-\gamma_i}$ for some exponents $\gamma_1,\gamma_2$
satisfying $\gamma_2>\gamma_1>0$.  \end{conjecture}

In the remainder of the paper we restrict our consideration to the
$d=1$ model.  In \csect{results} we state in \cthm{main} our main
result, the one-dimensional version of the conjecture, and in
\csect{measure} we describe the 1d limiting measure $\nu_\rho$.  The
proof of \cthm{main} is given in \csect{proof}.

\section{Statement of the result\label{results}}

The key to the rigorous establishment of a version of \cconj in
dimension $d=1$ is that there the existence and exact value of the
critical density are known---$\rho_c=1/2$---and that we also have a
complete description of the limiting measure $\nurho$ for $\rho\le1/2$
\cite{AGLS,GLS1,GLS2,ZC}.  (This measure was first identified and
discussed in \cite{ZC}, a reference which has just come to our
attention.  We regret that we did not properly credit this work in
earlier papers.)  We will discuss $\nurho$ for $\rho<1/2$ in
\csect{measure}; for the moment let us note that the measure at the
critical density is particularly simple:
 \be\label{muhalf}
\nu_{1/2}
   =\frac12(\delta_{\eta^*}+\delta_{\eta^\dagger}),
 \ee
 where $\eta^*$ and $\eta^\dagger$ are the two configurations in
$X=\{0,1\}^{\bbz}$ in which holes and particles strictly alternate.

 \Cconj concerns the asymptotic behavior of the quantity $V_\rho(L)$
with $L$ ``small,'' ``intermediate,'' or ``large,'' yet in each case
also sufficiently large.  To give a precise result in 1d we introduce
some notation for the description of this behavior.

\begin{definition}\label{asymp} Assume that $L$ is a positive integer
and $\delta$ a positive real number (in applications we will have
$\delta=\rho_c-\rho$), that $A(\delta,L)$ and $B(\delta,L)$ are real
valued functions whose asymptotic behavior in $L$ we wish to compare,
and that $L_1(\delta)$ and $L_2(\delta)$ are positive functions (these
play the role of setting the scales of the various regions).  Then we
write respectively
 \be\label{details}
A(\delta,L)\simeq B(\delta,L) \quad\text{for}\quad \begin{cases}
   \text{$L\ll L_1(\delta)$,}\\
   \text{$L_1(\delta)\ll L\ll L_2(\delta)$,}\\
   \text{$L\gg L_2(\delta)$,}\end{cases}
 \ee
 if for any $\epsilon>0$ there exists a $\delta_0>0$, together with a
(small) number $s>0$ and/or a (large) number $l>0$, such that for
$\delta<\delta_0$ we have, respectively,
 \be\label{details2}
1-\epsilon <\frac{A(\delta,L)}{B(\delta,L)}<1+\epsilon  \qquad \begin{cases}
   \text{for $L< sL_1(\delta)$,}\\
   \text{for $lL_1(\delta)< L< sL_2(\delta)$,}\\
   \text{for $L> lL_2(\delta)$.}\end{cases}
 \ee
   If $A$ and $B$ depend also on some additional parameter(s) $\alpha$
we say that \eqref{details} holds {\it uniformly} for $\alpha$ in some
(possibly $L$- or $\delta$-dependent) set if $s$ and/or $l$, as well
as $\delta_0$, may be chosen so that \eqref{details2} holds for all
such $\alpha$.  \end{definition}

With this notation established we may state our main result; we assume
that $0<\rho<1/2$, that $\delta=1/2-\rho$, and that $\nu_\rho$ is the
measure \eqref{nurho} of the 1d FEP.  In contrast to the $d\ge2$
behavior described in the \cconjn, for $d=1$ the variances $V_\rho(L)$
behave differently for $L$ odd and $L$ even.

\begin{theorem}\label{main}Let $\Nd(L)$ be the number of particles on
the sites $1,2,\ldots,L$, with distribution determined by $\nu_\rho$.   Then:

 \smallskip\noindent
 (a) For $L$ odd,
 \be\label{thma}
 V_\rho(L) =\Var(\Nd(L))
   \simeq\begin{cases}\frac14,&\text{for $L\ll\delta^{-2/3}$,}\\
  \frac23\sqrt{\frac2\pi}\delta L^{3/2},
          &\text{for $\delta^{-2/3}\ll L\ll\delta^{-2}$,}\\
  \frac14L,&\text{for $L\gg\delta^{-2}$.}\end{cases}
 \ee

 \smallskip\noindent
 (b) For $L$ even,
 \be\label{thmb}
V_\rho(L) =\Var(\Nd(L))
  \simeq\begin{cases}\frac23\sqrt{\frac2\pi}\delta L^{3/2},
          &\text{for $1\ll L\ll\delta^{-2}$,}\\
  \frac14L,&\text{for $L\gg\delta^{-2}$.}\end{cases}
 \ee
 \end{theorem} 

We actually have a stronger result for the asymptotics in the
large-$L$ region.  The ``right'' estimate for $\Var(\Nd(L))$ there is
$\rho(1-\rho)L$, as stated in \cconj and as we
discuss further in \crem{intuit} below, and this is true for all, not
just small, $\delta$:

\begin{corollary}[to the proof of \cthm{main}]
\label{stronger} For any $\epsilon>0$ there exists an $l>0$ such that
for any $\delta\in(0,1/2)$,
 \be\label{coreq}
 1-\epsilon < \frac{\Var(\Nd(L))}{\rho(1-\rho)L} < 1+\epsilon
   \qquad\text{for $L > l\delta^{-2}$.}
 \ee
 \end{corollary}

\noindent In fact, it follows from \eqref{forcor} below that
\eqref{coreq} holds with $l=(2\epsilon)^{-1}$.

 A comparison of \cthm{main}(a) with \cconj shows that for $L$ odd the
behavior of $V_\rho(L)$ in one dimension corresponds directly to the
conjectured behavior in higher dimension (but without the condition
$L\gg1$).  In particular it follows from \eqref{muhalf} that
 \be
 V_{1/2}(L)=\begin{cases}\frac14,& \text{if $L$ is odd,}\\
   0,& \text{if $L$ is even,}\label{v2l}\end{cases}
 \ee
 which implies that $\nu_{1/2}$ is hyperuniform and also explains the
$V_\rho(L)\simeq1/4$ behavior in \cthm{main} for small odd $L$.  The
variables introduced in the conjecture become $\lambda_1=0$,
$C_1=1/4$, $\lambda_2=3/2$, $C_2(\delta)=\sqrt{8/\pi}\,\delta/3$,
$\gamma_1=2/3$, and $\gamma_2=2$.  On the other hand, for even $L$ the
``small'' growth region is absent in one dimension: for small and
moderate values of $L$ the variances grow as $C_2(\delta)L^{3/2}$.
This odd/even distinction may be regarded as a legacy of \eqref{v2l}
when $\delta$ is perturbed away from 0.

\begin{remark}\label{intuit}
 \ccor{stronger} certainly implies that
$\lim_{L\to\infty}V_\rho(L)/L=\rho(1-\rho)$ for all $\rho<1/2$, and
this part of the result, although not the scale $\delta^{-2}$ at which
the limit is achieved, may be obtained by an elementary argument
\cite{GLS2}.  For with probability 1 each particle will move only a
finite distance during the evolution, so that for $L$ sufficiently
large $N(L)$ will, to high relative accuracy, be the same at the
end of the evolution as it was at the beginning, and  $\Var(N(L))$
will be the same as for the original Bernoulli measure.
\end{remark}

\begin{remark}\label{others} There are several one-dimensional models
with exclusion and facilitation, closely related to the FEP, for which
also $\rho_c=1/2$ and for which $\lim_{t\to\infty}\murho_t$ is for
$\rho\le1/2$ the measure $\nurho$ that we are considering here, and
hence for which the fluctuations $V_\rho(L)$ satisfy \cthm{main}.
In particular, this is true of the totally asymmetric, discrete-time
(parallel) updating in which all particles attempt to jump at the same
time, and only to the right \cite{GLS2}.  It is also true of an
asymmetric version of the continuous-time model of \csect{intro} in
which particles attempt to jump to the left or right at different rates
\cite{AGLS}.  \end{remark}

\section{The limiting measure $\nu_\rho$ for $\rho<1/2$\label{measure}}

A key ingredient for understanding the behavior described in
\cthm{main}, especially the behavior in the intermediate regime, is
the renewal structure of the stationary state $\nu_\rho$ defined in
\eqref{nurho}, $\rho<1/2$; we now describe this structure
\cite{GLS2,ZC}.  Since this state is frozen, the occurrence of
adjacent 1's has probability zero, so that the state is supported on
configurations of the form
 \begin{align}   \label{renew1}
  &\cdots1\ 0\ 1\ 0\ 1\ 0\ 1\ 0\ \hat0\ 1\ 0
   \ 1\ 0\cdots1\ 0\ \hat0\ \hat0\ \hat0\ 1\ 0\ 1\ 0\ 1\ 0
    \cdots1\ 0\ \hat0\cdots\\
  &\hskip50pt=\cdots0\ (1\ 0\ )^{X_{-1}}\ 0\ 
   (1\ 0)^{X_{0}}\ 0\ (1\ 0)^{X_{1}}\ 0
  \ (1\ 0)^{X_{2}}0\cdots.\label{renew2}
 \end{align}
 We focus particular attention on the $00$'s in \eqref{renew1}; 
the {\it second} $0$ of each such pair, marked as $\hat0$ in \eqref{renew1}
and corresponding to a 0 outside the parentheses in
\eqref{renew2}, will be called a {\it renewal event}.
(Note that adjacent renewal events correspond to a zero value for the
corresponding $X_i$.)  If we let $\hat\nu_\rho$ be
the measure $\nu_\rho$ conditioned on the occurrence of a renewal
event at the origin, then under $\hat\nu_\rho$ the $X_i$'s in
\eqref{renew2} are independent random variables that are identically
distributed, with the distribution described in \eqref{dist}.

The measure $\nu_\rho$ may now be described as the unique ergodic TI
measure such that conditioning on the occurrence of a renewal event at
the origin yields the measure $\hat\nu_\rho$ as just described.
Explicitly, if $\hat\nu_\rho^{(k)}$, $k\ge1$, denotes the restriction of
the measure $\hat\nu_\rho$ to configurations in which the first
renewal event to the right of the origin lies at site $k$, and
$\hat\nu_\rho^{(k,i)}$, $0\le i<k$, its translation by $i$ sites to
the left, then
$\nu_\rho=Z^{-1}\sum_{k\ge1}\sum_{i=1}^{k-1}\hat\nu_\rho^{(k,i)}$,
with $Z$ a normalizing factor.

\begin{remark}\label{dre} The density of the renewal events, i.e., the
probability  under $\nu_\rho$ of finding such an event at,
say, site $1$ (or equivalently of finding adjacent 0's at sites $0$
and $1$), is $1-2\rho=2\delta$, since under $\nu_\rho$ the probability
of adjacent 1's is zero.  \end{remark}

It is shown in \cite{GLS2,ZC} that the distribution of the $X_i$'s under
$\hat\nu_\rho$ is that of a random variable $\Xdh$ for which 
 \be\label{dist}\begin{aligned}
P(\Xdh=n)&=C_n\rho^n(1-\rho)^{n+1}\\
   &=\frac{1+2\delta}{2\cdot4^n}C_n(1-4\delta^2)^n,\end{aligned}
     \quad n=0,1,2,\ldots,
 \ee
 with $C_n$ the $n^{\rm th}$ Catalan number \cite{S}:
 \be
  C_n:=\frac1{n+1}\binom{2n}{n}
  =\frac{4^n}{n^{3/2}\sqrt\pi}\left(1+O\left(\frac1n\right)\right).
 \ee
  Here we have used Stirling's formula with error bounds. Thus for
$n\gg1$,
 \be\label{distasy}
P(\Xdh=n)\simeq\frac{1+2\delta}{2n^{3/2}\sqrt\pi}(1-4\delta^2)^n
  \simeq \frac{1}{2n^{3/2}\sqrt\pi}e^{-4\delta^2n},
 \ee
 where the first approximation holds for all $\delta$,
 $0\le\delta<1/2$, and the second for $\delta^2n\le1$.  (The
Catalan number $C_n$ arises here as the number of random walks of
length $2n$, with steps $\pm1$, which begin and end at the origin and
take only nonnegative values.)

It follows immediately from \eqref{muhalf} and \crem{dre} that
 \be\label{del0lim}
 \nu_{1/2}=\lim_{\delta\to0}\nu_\rho
 \ee
  in the sense of weak convergence, i.e., that
$\nu_{1/2}(A)=\lim_{\delta\to0}\nu_\rho(A)$ for every $A\subset X$
defined in terms of the configuration on a finite set of sites.  In
particular, the limiting measure contains no renewal events.  On the
other hand, the limit $\hat\nu_{1/2}=\lim_{\delta\to0}\hat\nu_\rho$ is
not so trivial: $\hat\nu_{1/2}$ is the probability distribution on
configurations of the form \eqref{renew2} for which there is a renewal
event at the origin and the i.i.d.~random variables $X_i$ have the
distribution of $\Xzh$:
 \be\label{limdist}
  P(\Xzh=n)=\lim_{\delta\to0}P(\Xdh=n)=\frac{C_n}{2\cdot4^n}
  \simeq\frac1{2\sqrt\pi n^{3/2}}.
  \ee
 In particular, there exist $c_1,c_2>0$ such that 
 \be\label{Pbounds}
 c_1n^{-3/2}<P(\Xzh=n)<c_2n^{-3/2},\qquad n\ge1.
 \ee
 Note that $\Xzh=\lim_{\delta\to0}\Xdh$ (limit in distribution); note
also that although $\hat\nurho$ for $\rho<1/2$ was obtained from
$\nurho$ by conditioning on a renewal event at the origin,
$\hat\nu_{1/2}$ cannot be so obtained from $\nu_{1/2}$.

$\Xzh$ has a $3/2$ power-law tail, and this 3/2 is, as we shall show,
the origin of the 3/2 in the $L^{3/2}$ behavior of the variance in the
intermediate regime.  (If 3/2 were replaced by $\gamma$, with
$1<\gamma\le2$, we would have had $L^\gamma$ behavior there
\cite{GLS4}).  Further, the fact that the exponential decay in
\eqref{distasy} becomes significant when $n$ is of order $\delta^{-2}$
is the origin of the fact that the transition to the large $L$ regime
occurs for $L$ of order $\delta^{-2}$.

 \medskip\noindent
 {\bf Notation:} Here, for the reader's convenience, we summarize our
notation, reviewing some that was introduced earlier and also defining
some new notation that will be used in the sequel.  We write
$\rho=1/2-\delta$, with $0\le\delta<1/2$, define
$\IL=\{1,2,\ldots,L\}$, and call the second of a pair of consecutive
empty sites a renewal event.  $\nu_\rho$ denotes the infinite-time
limit state \eqref{nurho} for the one-dimensional FEP at density
$\rho$, and $\hat\nu_\rho$ the state defined for $\delta>0$ by
conditioning $\nu_\rho$ on the occurrence of a renewal event at the
origin, and for $\delta=0$ as $\lim_{\delta\searrow0}\hat\nu_\rho$.  We
will use the following random variables:

 \smallskip\noindent
 \bll $\Nd(L)$, the number of particles in $\IL$ under $\nu_\rho$;

 \smallskip\noindent
 \bll $\Ndr(L)$ and $\Ndrh(L)$, the number of renewal events in $\IL$
under $\nu_\rho$ and $\hat\nu_\rho$, respectively (and, by convention,
$\Ndrh(0)=0$);

 \smallskip\noindent
\bll $\Xdh$, a random variable with the distribution under $\hat\nu_\rho$
of the $X_i$ in \eqref{renew2};

 \smallskip\noindent
\bll $\Ydh=2\Xdh+1$, the distance between renewal events under
$\hat\nu_\rho$;

\bll $\Fd$, the location of the first renewal event to
the right of the origin under $\nu_\rho$.

We make a further remark and a general notational
convention.  First, the translation invariance of $\nurho$ implies
that the number of renewal events in any interval of $L$ sites,
conditioned on the occurrence of a renewal event at the immediately
preceding site, has the same distribution as $\Ndrh(L)$.  Second, we
will let $P$ denote probabilities for random variables such as those
above, whether defined using $\nu_\rho$ or $\hat\nu_\rho$; thus, for
example, $P(\Ndr(L)=n)=\nu_\rho(\Ndr(L)=n)$ while
$P(\Ndrh(L)=n)=\hat\nu_\rho(\Ndrh(L)=n)$.  

 Note that from \crem{dre},
 \be\label{ENr}
  E\bigl(\Ndr(L)\bigr)=2\delta L.
 \ee
Note also that while, from \eqref{del0lim},
$\Nzr(L)=\lim_{\delta\to0}\Ndr(L)$ (limit in distribution) is the zero
random variable, since in $\nu_{1/2}$ there are no renewal events,
$\Nzrh(L)=\lim_{\delta\to0}\Ndrh(L)$ is non-trivial.  We stress that
the key to the $\delta\searrow0$ asymptotics described in \cthm{main}
lies not in $\nu_{1/2}$ but in $\hat\nu_{1/2}$.

\section{Proof of \cthm{main}\label{proof}}

The proof of \cthm{main} is broken into nine steps, as follows:

\begin{enumerate}[{\bf Step 1:}]

\item Express $\Nd(L)$ in terms of $\Ndr(L)$.

\item Express the distribution and second moment of
$\Ndr(L)$ in terms of those of $\Ndrh(L)$.

\item Approximate the expressions found in Step~2 by replacing
  $\Ndrh(L)$ by $\Nzrh(L)$.

\item The expressions obtained in Step~3 involve the quantity
$p_\delta$ defined in \eqref{S21}. In Step~4 we express this quantity
in terms of $\Ydh$, then replace the occurrences of $\Ydh$ by $\Yzh$.

\item Obtain the large-$L$ asymptotics of the distribution of $\Nzrh(L)$ and of
its second moment, and insert these into the expressions found in
Step~4.

\item Use the asymptotics of the
distribution of $\Yzh$ to further approximate the expressions found in
Step~5.

\item Obtain from the expressions found in Step~6 the asymptotics of
  $\Var(\Ndr(L))$.

\item Use Step~1 to obtain the results of \cthm{main} for
$L\ll\delta^{-2}$ from the expression found in Step~7 for
$\Var(\Ndr(L))$.

\item Use some facts about the truncated two-point correlation
  function for $\nu_\rho$ to obtain the results of \cthm{main} for
$L\gg\delta^{-2}$.

\end{enumerate} 

We now consider these steps in order.

 \medskip\noindent
     {\bf Step 1:} In this step we again use the notation introduced
in \eqref{renew1}, so that the values which may be taken by a
configuration $\eta_i$ are $\hat0$, $0$, and $1$, where $\hat0$
denotes a renewal event, $0$ an empty site preceded by a $1$, and $1$
an occupied site.  Now we observe that $L-\Ndr(L)$ is odd if and only
if the pair $(\eta_1,\eta_L)$ has value $(0,0)$, $(0,\hat0)$,
$(\hat0,1)$, or $(1,1)$; moreover,
 \be\label{Nsigg}
\Nd(L)=\frac12\bigl[L-(\Ndr(L)+\sd(L))\bigr],
 \ee
 where
 \be\label{sddef}
\sd(L)=\begin{cases}0,&\text{if $L-\Ndr(L)$ is even,}\\
1,&\text{if $L-\Ndr(L)$ is odd and
  $(\eta_1,\eta_L)$ is $(0,0)$ or $(0,\hat0)$,}\\
-1,&\text{if $L-\Ndr(L)$ is odd and 
  $(\eta_1,\eta_L)$ is $(\hat0,1)$ or $(1,1)$.}
\end{cases}
 \ee
 One checks this by induction on $\Ndr(L)$; the case $\Ndr(L)=0$ is
easy.  For the induction step one
passes from a configuration $\eta$ to another $\eta'$ by removing a
$\hat0$ from some site $i$ with $1\le i\le L$ and setting
$\eta'_j=\eta_j$ if $j<i$ and $\eta'_j=\eta_{j+1}$ if $j\ge i$; one
then applies the induction assumption to $\eta'$ on $J_{L-1}$, noting
that then $L$ and $\Ndr(L)$ both decrease by 1, and observing that
$(\eta_1',\eta_{L-1}')\ne(\eta_1,\eta_L)$ only if $i=1$ and
$\eta_1'=1$ or $i=L$ and $\eta'_{L-1}=0$, and that in each of these
cases $\sigma$ is unchanged and \eqref{Nsigg} remains valid.

From \eqref{Nsigg} we have that 
 \be\label{varnd}
\Var(\Nd(L)=\frac14\Var\bigl(\Ndr(L)+\sd(L)\bigr).
 \ee
 To simplify this expression further we note that, writing $E$ for
 expectation, we have
 \be\label{Econd}
E\bigl(\sd(L)\,\big|\,\Ndr(L)=n\bigr)=0, \quad\text{for any $n\ge0$,}
 \ee
 as we will argue shortly.  But then
 \be\label{Esig}
E(\sd(L))=0
 \ee
and 
 \be
 E(\Ndr(L)\sd(L))=\Cov(\Ndr(L),\sd(L))=0,
 \ee
 so that
 \be
\Var\bigl(\Ndr(L)+\sd(L)\bigr) 
   =\Var\bigl(\Ndr(L)\bigr)+\Var\bigl(\sd(L)\bigr).
 \ee
 Thus from \eqref{varnd} and \eqref{Esig},
 \begin{align}\nonumber
\Var(\Nd(L))&=\frac14\big(\Var(\Ndr(L))+\Var(\sd(L))\bigr)\\
 \label{varnd2}
  &=\frac14\big(\Var(\Ndr(L))+P(\sd(L)\ne0)\bigr).
  \end{align}

To verify \eqref{Econd} we first note that, from \eqref{nurho} and the
reflection invariance of the Bernoulli measure and of the dynamics,
$\nu_\rho$ is invariant under reflection about any (integer or
half-integer) point.  \eqref{Econd} then follows from the observation
that reflection about $L/2$, the midpoint of the
interval $\{0,1,\ldots,L\}$ ({\it not} of $J_L$), leaves $\Ndr(L)$
unchanged and, when $L-\Ndr(L)$ is odd, changes the sign of $\sd(L)$,
as one sees by checking separately for the four possible values of
$(\eta_1,\eta_L)$ which can then occur.

 \medskip\noindent
 {\bf Step 2:} With $\Fd$ as defined towards the end of
\csect{measure}, let 
 \be\label{S21}
\pd(l)=P(\Fd=l).
 \ee
   The critical
observation for Step 2 is that the number of renewal events in $\IL$,
if there are any, is one more than the number of such events to the
right of $\Fd$.  Thus for $n\ge1$,
 \begin{align}\nonumber
P\bigl(\Ndr(L)=n\bigr)
  &=\sum_{l=1}^L\pd(l)P\bigl(\Ndr(L)=n\,\big|\,\Fd=l)\bigr)\\
   &=\sum_{l=1}^L\pd(l)P\bigl(\Ndrh(L-l)=n-1\bigr),\label{useF}
 \end{align}
 and so
 \be\label{S23}
  E\bigl(\Ndr(L)^2\bigr)=\sum_{l=1}^L\pd(l)E\bigl((\Ndrh(L-l)+1)^2\bigr).
 \ee

 \medskip\noindent
 {\bf Step 3:} As indicated earlier, the next step is to control the
approximation arising from the replacement of $\Ndrh$ by $\Nzrh$ in
\eqref{useF} and \eqref{S23}.  Specifically, we will show that for
$n\ge1$,
 \be\label{S26}
P\bigl(\Ndr(L)=n\bigr)
   \simeq\sum_{l=1}^L\pd(l)P\bigl(\Nzrh(L-l)=n-1\bigr)\
   \quad\text{for $L\ll\delta^{-2}$,}
 \ee
 uniformly in $n\le k\sqrt L$ for $k$ any fixed positive integer (see
\cdef{asymp}), and also that
 \be\label{S28}
  E\bigl(\Ndr(L)^2\bigr)\simeq\sum_{l=1}^L\pd(l)E\bigl((\Nzrh(L-l)+1)^2\bigr)
       \quad\text{for $L\ll\delta^{-2}$.}
 \ee
 Note that the right hand sides of \eqref{S26} and \eqref{S28} both
mix quantities defined for $\delta>0$ with those defined for
$\delta=0$ ($p_\delta$ and $\Nzrh$ respectively).  These equations are
more delicate than they may appear because they demand that we control
the errors in these approximations by requiring merely that, for small
$\delta$, the quantity $L\delta^2$ be sufficiently small
regardless of the size of $L$ itself.

Let us fix $\Lt< L$ ($\Lt$ plays the role of $L-l$ in \eqref{S26}
and \eqref{S28}) and let $\etal$ and $\etaldh$ denote respectively a
fixed and a random configuration on $\ILt$, such that (i)~$\etal$
contains $m=n-1$ renewal events, with the convention that if $\etal_1=0$
then this is a renewal event, and (ii)~$\etaldh$ is the restriction to
$\ILt$ of a configuration distributed according to $\hat\nu_\rho$.
Then for $\delta\ge0$,
 \be\label{Nsum}
 P\bigl(\Ndrh(\Lt)=n-1\bigr)=\sum_{\etal}P(\etaldh=\etal).
 \ee
 \eqref{S26} will now follow from \eqref{Nsum} once we show that
 for all $\etal$, uniformly in $m\le k\sqrt L$, 
 \be\label{etaeta}
 P(\etaldh=\etal)\simeq P(\etalzh=\etal) 
       \quad\text{for $L\ll\delta^{-2}$.}
 \ee

 To verify \eqref{etaeta} we let $2x_i+1$, $i=1,\ldots,m$, be the distances
between the renewal events in $\etal$, with $2x_1+1$ the distance from
the origin to the first renewal event.  Further, we define
 \be\label{qdeldef}
\qd(L)=P\bigl(\Ydh>L\bigr)
=P\left(\Xdh>\left\lfloor\frac{L-1}2\right\rfloor\right)  
 \ee
  with $\Ydh=2\Xdh+1$ as defined at the end of \csect{measure}, and
 \be\label{rdeldef}
 \rd(L)=\frac{\qd(L)}
   {\qz(L)}.
 \ee
 Now  from \eqref{dist} and \eqref{limdist},
 \be\label{newdist}
P(\Xdh=l)=(1+2\delta)(1-4\delta^2)^lP(\Xzh=l).
 \ee
 Then, with $L'$ the distance from the last renewal event in $\etal$
 to the right boundary $\Lt$ of $J_{\Lt}$, we have 
 \begin{align}\nonumber
P(\etaldh=\etal)
   &=\prod\nolimits_{i=1}^mP(\Xdh=x_i)\ \qd(L')\\
   &=(1+2\delta)^m(1-4\delta^2)^{\sum_{i=1}^mx_i}
   \rd(L')P(\etalzh=\etal).\label{S30}
 \end{align}

We complete the argument for \eqref{etaeta} by showing that each of
the first three factors on the right hand side of \eqref{S30} is
asymptotic to 1 for $L\ll\delta^{-2}$, uniformly in $n\le k\sqrt{L}$
for any fixed positive $k$.  First, under this condition,
 \be\label{etf1}
1\le(1+2\delta)^m \le(1+2\delta)^{k\sqrt L}\le e^{2k\delta\sqrt L}.
 \ee
Next, from $1-x\ge e^{-2x}$ for $0\le x\le1/2$ we have for any $M\ge0$,
 \be\label{etf2}
 1\ge(1-4\delta^2)^M\ge e^{-8M\delta^2} \ge1-8M\delta^2
  \quad\text{for $0\le\delta\le\frac1{2\sqrt2}$},
 \ee
 and with $M=\sum_{i=1}^mx_i\le L/2$ this gives the desired
asymptotics of the second factor.  We will finally show that 
 \be\label{etf3}
  \rd(L)\simeq1 \quad\text{for $L\ll\delta^{-2}$};
 \ee
 this, with $L$ replaced by $L'$, gives the needed control of the right
 hand side of \eqref{S30}.

 We first observe that \eqref{Pbounds} implies that there exists a
$c>0$, independent of $L$, such that if for $\epsilon>0$ we set
$K= c\epsilon^{-2}L$ then the tail of the series for $q_0(L)$ beyond
$K$ is relatively small, i.e., $q_0(K)<\epsilon q_0(L)$, which implies
 \be\label{new1}
 q_0(L) < \frac1{1-\epsilon}\sum_{L<l\le K} P(\Yzh=l).
 \ee
 Further, for  $L\delta^2<\epsilon^3/(8c)$ it follows from
\eqref{etf2} that
 \be\label{new2}
 \sum_{L<l\le K}P(\Ydh=l) > (1-\epsilon)\sum_{L<l\le K}P(\Yzh=l).
 \ee
 Then from \eqref{dist}, \eqref{new1}, and \eqref{new2},
 \be\label{new3}
 1+2\delta\ge\frac{q_\delta(L)}{q_0(L)}
  >(1-\epsilon)\frac{\sum_{L<l\le K}P(\Ydh=l)}{\sum_{L<l\le
      K}P(\Yzh=l)}
 > (1-\epsilon)^2.
 \ee
 This completes the verification of \eqref{etf3}.  \eqref{etaeta} then follows from
 \eqref{S30}--\eqref{etf3}.

We remark that in fact the first inequality in \eqref{new3} can be
strengthened:
 \be\label{new4}
 \frac{q_\delta(L)}{q_0(L)}\le1\qquad \text{for all $L\ge1$.}
 \ee
 For from \eqref{dist} the difference
$q_0(L)-q_\delta(L)$ is increasing for $L\le L^*$ and decreasing
for $L>L^*$, where $L^*=-\log(1+2\delta)/\log(1-4\delta^2)$, and
vanishes at $L=0$ and as $L\to\infty$.

 We now turn to \eqref{S28}.  For any positive $k$ we write
$\Ndrh(L)=\Ndlh{k\sqrt L}(L)+\Ndgh{k\sqrt L}(L)$, where
 \be
 \Ndlh{x}(L):=\Ndrh(L) I_{\bigl\{\Ndrh(L)\le x\bigr\}}
 \ee
 and
 \be
  \Ndgh{x}(L):=\Ndrh(L) I_{\bigl\{\Ndrh(L)> x\bigr\}},
 \ee
 with $I_{\{\cdot\}}$ denoting the indicator function of the set
$\{\cdot\}$.  Then  from \eqref{etaeta},
 \be\label{rhsl}
   E(\Ndlh{k\sqrt L}(L)^2)\simeq E(\Nzlh{k\sqrt L}(L)^2) 
  \quad\text{for $L\ll\delta^{-2}$.}
 \ee
  \eqref{S28} will follow easily once we strengthen \eqref{rhsl} to
 \be\label{goal}
   E(\Ndrh(L)^2)\simeq E(\Nzrh(L)^2) \quad\text{for $L\ll\delta^{-2}$.}
 \ee

There are two crucial facts for obtaining \eqref{goal} from
\eqref{rhsl}, via the approximations expressed in \eqref{twoeq} below.  The
first, to be proved shortly, is that for any $\epsilon>0$ there exists
a $k>0$ such that
 \be\label{S38}
 E\bigl(\Ndgh{k\sqrt L}(L)^2\bigr)\le\epsilon L \quad
  \text{for $L\ll\delta^{-2}$}
 \ee
 (where $L\ll\delta^{-2}$ holds for all $L$ if $\delta=0$). The second,
to be proved in Step~5, is that $\lim_{L\to\infty}E(\Nzrh(L)^2)/L=1$, so
that for some constant $C>0$,
 \be\label{get282}
  E(\Nzrh(L)^2)\ge CL
  \quad\text{for all $L\ge1$.}
 \ee
 \eqref{goal} follows from \eqref{rhsl}, \eqref{S38} (for both
$\delta=0$ and $\delta>0$) and \eqref{get282}.

  To see this, fix $\epsilon>0$ and take $k$ so that \eqref{S38}
holds.  Then from \eqref{S38} for $\delta=0$ and \eqref{get282} we
have (for all $L$)
 \be\label{get283}
\left|\frac{E\bigl(\Nzlh{k\sqrt L}(L)^2\bigr)}
  {E\bigl(\Nzrh(L)^2\bigr)}-1\right|
 =\frac{E\bigl(\Nzgh{k\sqrt L}(L)^2\bigr)}
  {E\bigl(\Nzrh(L)^2\bigr)}<\frac\epsilon C.
 \ee
  Now by \eqref{rhsl} we may take $L\delta^2$ so small that 
 \be\label{get2}
\left|\frac{E\bigl(\Ndlh{k\sqrt L}(L)^2\bigr)}
  {E\bigl(\Nzlh{k\sqrt L}(L)^2\bigr)}-1\right| <\epsilon. \ee
From this,
 \begin{align}\nonumber
 E\bigl(\Ndlh{k\sqrt L}(L)^2\bigr)
  &>E\bigl(\Nzlh{k\sqrt L}(L)^2\bigr)(1-\epsilon)\\
  &=\label{get3}\big[E\bigl(\Nzrh(L)^2\bigr)
     -E\bigl(\Nzgh{k\sqrt L}(L)^2\bigr)\big](1-\epsilon)\\
  &>(C-\epsilon)(1-\epsilon)L,\nonumber
 \end{align}
 and then  we have from  \eqref{S38} (possibly with a further
 restriction on $L\delta^2$) that
 \be\label{get4}
\left|\frac{E\bigl(\Ndlh{k\sqrt L}(L)^2\bigr)}
  {E\bigl(\Ndrh(L)^2\bigr)}-1\right|
 \le\frac{E\bigl(\Ndgh{k\sqrt L}(L)^2\bigr)}
  {E\bigl(\Ndlh{k\sqrt L}(L)^2\bigr)}
  <\frac\epsilon{(C-\epsilon)(1-\epsilon)}.
 \ee
 Since $\epsilon$ here is arbitrary, \eqref{get283} and \eqref{get4}
imply respectively that
 \be\label{twoeq}
 E\bigl(\Nzrh(L)^2\bigr)\simeq E\bigl(\Nzlh{k\sqrt L}(L)^2\bigr)
\quad \text{and}\quad
 E\bigl(\Ndrh(L)^2\bigr)\simeq E\bigl(\Ndlh{k\sqrt L}(L)^2\bigr)
 \ee
 for $L\ll\delta^{-2}$, and with \eqref{rhsl}, \eqref{goal} follows.
 
To conclude Step~3 we must establish \eqref{S38}.  To do so we first
note that for any integer $n\ge1$,
 \be\label{S41}
  P(\Ndrh(L)\ge n)\le P\left(\Ydh\le L\right)^n
  = (1-\qd(L))^n,
 \ee
  since the distances between the successive renewal events in $\IL$
(including the distance of the first such event from the origin),
which are independent, must each be no greater than $L$.

In the remainder of this section we write $q=\qd(L)$.  For
any integer-valued random variable $N$, and any integer $n_c\ge1$, we
have (as a consequence of summation by parts) that
 \be\label{S42}
E(N^2I_{\{N\ge n_c\}})
  =n_c^2P(N\ge n_c)+\sum_{n=n_c+1}^\infty(2n-1)P(N\ge n),
 \ee
 provided that $n^2P(N\ge n)\to0$ as $n\to\infty$.  Thus from
\eqref{S41},
 \begin{align}\nonumber
 E(\Ndgh{n_c}(L)^2)
 &\le n_c^2(1-q)^{n_c}+\sum_{n=n_c+1}^\infty(2n-1)(1-q)^n\\
  &\le \left(n_c^2+\frac{2n_c}q+\frac2{q^2}\right)(1-q)^{n_c}\nonumber\\
  &\le2\left(n_c+\frac1q\right)^2e^{-qn_c}.\label{S43}
 \end{align}
 If $\delta>0$ then it follows from \eqref{distasy} and
\eqref{qdeldef} that for $L\delta^2$ sufficiently small there is an
$A>0$ such that
 \be\label{S44}
 q\ge\frac{A}{\sqrt{L}}.
 \ee
 Moreover, \eqref{limdist} implies the same conclusion, for any $L$,
when $\delta=0$.  If now for any $k>0$ we set 
$n_c=\left\lfloor k\sqrt L\right\rfloor$ then \eqref{S43} and
\eqref{S44} yield
 \be
 E(\Ndgh{k\sqrt L}(L)^2)=E(\Ndgh{n_c}(L)^2)
    \le\left(k+\frac{1}{A}\right)^2e^Ae^{-Ak}\,L,
 \ee
 and \eqref{S38} will hold for sufficiently large $k$.

 \medskip\noindent
 {\bf Step 4:} Equation \eqref{S28}, the starting point for our future
investigations, involves $p_\delta(l)=P(\Fd=l)$, the probability under
$\nu_\rho$ that the first renewal event to the right of the origin
occurs at site $l>0$.  This happens precisely when there is a renewal
event at some site $-l'\le0$, an event with probability $2\delta$
(see \crem{dre}), and the next renewal event to its right
is at $l$, so that
 \be
p_\delta(l)=2\delta\sum_{l'\ge0}P(\Ydh=l'+l)
  =2\delta P(\Ydh\ge l).
 \ee
 (Note that $p_\delta(l)=2\delta\qd(l-1)$.)  Now \eqref{etf3} yields
 \be\label{S46}
p_\delta(l)\simeq2\delta P(\Yzh\ge l)  
  \quad\text{for $l\ll\delta^{-2}$,}
 \ee
 and thus from \eqref{S28} we have that 
 \be\label{S48}
    E\bigl(\Ndr(L)^2\bigr)
  \simeq2\delta\sum_{l=1}^L P(\Yzh\ge l) E\bigl((\Nzrh(L-l)+1)^2\bigr)
       \quad\text{for $L\ll\delta^{-2}$.}
 \ee

 \medskip\noindent
 {\bf Step 5:} In this step we obtain the
large-$L$ asymptotics of $E\bigl(\Nzrh(L)^2\bigr)$: 
 \be\label{S51}
E\bigl(\Nzrh(L)^2\bigr)\simeq L, \quad\text{for $L\gg1$.}
 \ee
 This formula may be obtained from \cite{GL}, but we give a
self-contained proof, arguing from the detailed form of the
distribution of the renewal random variable $\Yzh=2\Xzh+1$.  Recall
(see \eqref{dist}) that $P(\Xzh=n)=C_n2^{-(2n+1)}$; the Catalan number
$C_n$ counts the number of paths between time 0 and time $2n$ of a
random walk which starts and ends at the origin while never taking any
positive value.  Thus $\Yzh$ has the same distribution as the time of
first arrival at site 1 of a simple symmetric random walk $W_l$,
$l=0,1,\ldots$, which starts at the origin.  As a consequence,
$\Nzrh(L)$ has the same distribution as the maximum value $M(L)$ of
$W_l$ over the interval $[0,L]$.  From \cite{F}, Section III.7,
Theorem 1 we then have
 \be\label{frfel}
 P(\Nzrh(L)=n)
  =\begin{cases}
  P(W_L=n),&\text{if $L-n$ is even,}\\
  P(W_L=n+1),&\text{if $L-n$ is odd.}
    \end{cases}
 \ee
 An easy calculation from \eqref{frfel} gives, for $L$ odd,
 \be\label{easy}
E\bigl(\Nzrh(L)^2\bigr)
   =E\bigl(W_L^2\bigr)-E\bigl(|W_L|\bigr) 
   +\frac12\bigl(1-P(W_L=0)\bigr),
 \ee
 and this yields \eqref{S51}, since $E\bigl(W_L^2\bigr)=L$ and
$E\bigl(|W_L|\bigr)\le\sqrt L$, by the Cauchy-Schwarz inequality.  

\begin{remark}\label{step5} (a) From \eqref{frfel} one can show easily
that $\Nzrh(L)/\sqrt L$ converges in distribution, as $L\to\infty$,
to $|Z|$, with $Z$ a standard normal random variable.

 \smallskip\noindent
 (b) In \cite{GLS2,ZC} a random walk representation of particle
configurations (there called a height function or height process) was
used to obtain the distribution \eqref{dist}.  The Catalan numbers
play the same role in this derivation that they do above.

\end{remark}

 \medskip\noindent
 {\bf Step 6:} Equation \eqref{S48} provides the leading order
small-$\delta$ approximation to $E\bigl(\Ndr(L)^2\bigr)$, valid for
$L\ll\delta^{-2}$. In this step we use \eqref{S51} and \eqref{qasym}
below to approximate this moment to leading order in $L$, as well.

First, note that from \eqref{distasy} we have at once that
 \be\label{qasym}
 P\bigl(\Yzh\ge l\bigr)=P\Bigl(\Xzh\ge\frac{l-1}2\Bigr) 
  \simeq\sqrt{\frac2{\pi l}} \qquad\text{for $l\gg1$.}
 \ee
 Substituting \eqref{S51} and \eqref{qasym} into \eqref{S48} yields, at
least formally,
 \be\label{S53}
E\bigl(\Ndr(L)^2\bigr)\simeq2\sqrt{\frac2\pi}\delta\sum_{l=1}^L
\frac{L-l}{\sqrt l}, \quad\text{for $1\ll L\ll\delta^{-2}$.}
 \ee
 We will justify \eqref{S53} shortly, but for the moment only note
that the restriction $L\gg1$, not present in \eqref{S48}, arises from
\eqref{S51} and \eqref{qasym}.  From \eqref{S53} we have that
  for $1\ll L\ll\delta^{-2}$,
 \begin{align}\label{S70}
E\bigl(\Ndr(L)^2\bigr)\nonumber
  &\simeq 2\sqrt{\frac2\pi}\delta\int_1^L\frac{L-x}{\sqrt x}\,dx\\
   &= 2\sqrt{\frac2\pi}\delta L^{3/2}\int_{1/L}^1\frac{1-y}{\sqrt y}\,dy
  \simeq \frac83\sqrt{\frac2\pi}\delta L^{3/2}.
 \end{align}
 This is our final approximation for $E\bigl(\Ndr(L)^2\bigr)$.

We now return to \eqref{S53}.  We are justified (when $l$ is not too
large) in replacing $E\bigl((\Nzrh(L-l)+1)^2\bigr)$ by
$E\bigl((\Nzrh(L-l))^2\bigr)$ in passing from \eqref{S48} to
\eqref{S53} since, from \eqref{S51} and the Cauchy-Schwarz inequality,
$E\bigl(\Nzrh(L)\bigr)\ll E\bigl(\Nzrh(L)^2\bigr)$ for $L\gg1$.
Moreover, the substitutions in \eqref{S48} suggested by \eqref{S51}
and \eqref{qasym} are valid provided that $l$ and $L-l$ are
sufficiently large.  More precisely,  using \eqref{S51} and
\eqref{qasym}, we can conclude that there exists an integer $l_*$,
which does not depend on $L$, such that
 \be
 2\sqrt{\frac2{\pi}}\delta\sum_{l=l_*}^{L-l_*}\frac{L-l}{\sqrt l} 
 \simeq 2\delta\sum_{l=l_*}^{L-l_*} P(\Yzh\ge l) 
 E\bigl((\Nzrh(L-l)+1)^2\bigr).
 \ee
 But note that the sums in \eqref{S48} and \eqref{S53} over
$1\le l<l_*$ and $L-l_*<l\le L$ are $O(L)$ and $o(1)$,
respectively, as $L\to\infty$, while the sums over
$l_*\le l\le L-l_*$ is of order $L^{3/2}$ (see
\eqref{S70}).  This completes the justification.

 \medskip\noindent
  {\bf Step 7:} From \eqref{ENr} and \eqref{S70},
 \be
 \frac{E\bigl(\Ndr(L)\bigr)^2}{E\bigl(\Ndr(L)^2\bigr)}
 \simeq\frac32\sqrt{\frac\pi2}\delta L^{1/2}\ll1
  \quad\text{for $1\ll L\ll\delta^{-2}$.}
 \ee
 Thus, again from \eqref{S70},
 \be\label{S71}
 \Var\bigl(\Ndr(L)\bigr) \simeq E\bigl(\Ndr(L)^2\bigr)
  \simeq \frac83\sqrt{\frac2\pi}\delta L^{3/2},
   \quad\text{for $1\ll L\ll\delta^{-2}$.}
 \ee

 \medskip\noindent
 {\bf Step 8:} We can now combine the results of Step~1 and Step~7 to
obtain the parts of \cthm{main} which concern $L\ll\delta^{-2}$.
For from \eqref{varnd2} and \eqref{S71} we have that
 \be\label{S75}
 \Var(\Nd(L))\simeq\frac23\sqrt{\frac2\pi}\delta L^{3/2}
   +\frac14 P(\sd(L)\ne0)\quad\text{for $1\ll L\ll\delta^{-2}$,}
 \ee
  with $\sd(L)$ defined in \eqref{sddef}.

If $L$ is even then, from \eqref{ENr},
 \be\label{S76}\begin{aligned}
P(\sd(L)\ne0)&=P(\Ndr(L)\text{ is odd})\\
 &\le P(\Ndr(L)>0) \le E(\Ndr(L)))=2\delta L.
 \end{aligned}\ee
 Since $L\ll L^{3/2}$ for $L\gg1$, the $L\ll\delta^{-2}$ part of
\cthm{main}(b) follows.

 On the other hand, if $L$ is odd then, again from \eqref{ENr},
 \be\label{S77}\begin{aligned}
P(\sd(L)\ne0)&=P(\Ndr(L)\text{ is even})\\
 &\ge P(\Ndr(L)=0) \ge 1-2\delta L\simeq1\quad\text{for $L\ll\delta^{-1}$.}
 \end{aligned}\ee
 Since $\delta L^{3/2}\gg1$ for $L\gg\delta^{-2/3}$ and
$\delta L^{3/2}\ll1$ for $L\ll\delta^{-2/3}$, we obtain from
\eqref{S77} and \eqref{S75} the conclusions of \cthm{main}(a) for
$1\ll L\ll\delta^{-2}$.  To remove the restriction that $L\gg1$ we
note that
 \be
 \Var\bigl(\Ndr(L)\bigr)\le E(\Ndr(L)^2)\le E(\Ndr(L')^2)
 \ee
 for $L\le L'$. Choosing $L'$ such that also
$1\ll L'\ll\delta^{-2/3}$, we see using \eqref{varnd2}, \eqref{S70},
and \eqref{S77} that $\Var(\Nd(L))\simeq\frac14$ for
$1\le L\ll\delta^{-2/3}$.

\begin{remark}\label{Step8}Concerning the estimate
$P(\Ndr(L)>0)\le2\delta L$ used in \eqref{S76} and \eqref{S77}, note
that it follows from \eqref{S21}, \eqref{S46}, and \eqref{qasym} that
in fact for $L\gg1$,
 \be
  P(\Ndr(L)>0)=P(\Fd\le L)\simeq2\delta\int_0^L\sqrt{\frac2{\pi l}}\,dl
  \simeq 4\sqrt{\frac2\pi}\delta\sqrt L.
 \ee
 \end{remark}

 \medskip\noindent
 {\bf Step 9:} We now turn to the $L\gg\delta^{-2}$ part of
\cthm{main} and to the related \ccor{stronger} (recall also
\crem{intuit}).  
We write
 \begin{align}
\Var(\Nd(L))&=\sum_{1\le k\le L}\Var\nolimits_{\nu_\rho}(\eta_k)
  +2\sum_{1\le k<i\le L}\Cov\nolimits_{\nu_\rho}(\eta_k,\eta_i)\nonumber\\
 &=\rho(1-\rho)L+2\sum_{k=1}^{L-1}\sum_{j=1}^k\gtr(j),\label{S83}
 \end{align}
 with $\gtr$ the truncated two-point correlation function for the TI
state $\nu_\rho$:
 \be
\gtr(k)=E(\eta_j\eta_{j+k})-\rho^2.
 \ee
 It is shown in \cite{GLS2} that for all $n\ge0$,
 \be
 \gtr(2n+1)+\gtr(2n+2)=0,
 \ee
  so that \eqref{S83} becomes 
 \be\label{S85}
\Var(\Nd(L))=\rho(1-\rho)L
  +2\sum_{\latop{k=1}{k\text{ odd}}}^{L-1}\gtr(k).
 \ee
  We will show shortly that for  $j\ge0$, 
 \be\label{S88}
 \bigl|\gtr(2j+1)\bigr|\le\rho^2 (1-4\delta^2)^j.
 \ee
 Then from \eqref{S85},
 \be\label{forcor}
 \left|\frac{\Var(\Nd(L))}{\rho(1-\rho)L}-1\right| 
 \le\frac{2\rho}{(1-\rho)L}\sum_{j=0}^\infty(1-4\delta^2)^j
    <\frac1{2\delta^2L},
 \ee
 and this verifies the result stated in \ccor{stronger}.  The
$L\gg\delta^{-2}$ cases of \cthm{main} then follow by taking the $\delta_0$
 of \cdef{asymp} sufficiently small.

  Consider now \eqref{S88}.   The generating function of the $\gtr$
is computed in \cite{GLS2}:
 \be\label{gfct}
 G_\rho(z):=\sum_{k=1}^\infty\gtr(k)z^k
 =\frac{z\bigl(\sqrt{1-z^2(1-4\delta^2)}-2\delta\bigr)^2}
   {4(z-1)(z+1)^2}.
 \ee
 The numerator in \eqref{gfct} has double zeros (on its first sheet)
at $z=\pm1$, so that $G_\rho$ is analytic at these points with
singularities at $z=\pm z_*$, where $z_*:=(1-4\delta^2)^{-1/2}$.
Expressing $\gtr(k)$ via Cauchy's formula as an integral over a small
circle around the origin, distorting this contour to obtain the sum of
integrals of the discontinuity of $G_\rho$ across cuts on the real
axis from $z_*$ to $\infty$ and from $-z_*$ to $-\infty$, and making
the change of variable $z\to-z$ in the second of these integrals, we
obtain a representation of $\gtr(k)$ which for $k$ odd is
 \be
\gtr(k)=-\frac{2\delta}{\pi z_*}
  \int_{z_*}^\infty\frac{\sqrt{z^2-z_*^2}}{z^k(z^2-1)^2}\,dz,\qquad
  \text{$k$ odd.}
 \ee
 This implies that for $k$ odd, $|\gtr(k)|\le z_*^{-(k-1)}|\gtr(1)|$,
and since $\gtr(1)=-\rho^2$, \eqref{S88} follows.

This completes the proof of \cthm{main}.

\section{Concluding remarks\label{conclude}}

We note that on $\bbz$ the approach to hyperuniformity as
$\rho\searrow\rho_c$ is different from that for $\rho\nearrow\rho_c$;
in the former case the unique stationary measure for the FEP on $\bbz$
is known and \cite{LY,GLS4}
 \be
\lim_{L\to\infty}\frac1LV_\rho(L)=\rho(1-\rho)(2\rho-1) 
  \quad\text{for $\rho>1/2$.}
 \ee
Thus from \cthm{main} and \eqref{v2l},
$\lim_{L\to\infty}L^{-1}V_\rho(L)$ is continuous in $\rho$ from above,
but not from below, at $\rho=1/2$.  We expect similar behavior on
$\bbz^d$.

 The case of the FEP on a {\it ladder}, a system consisting of two
(infinite) rows of sites, was studied numerically in \cite{O,I}.  Here
again, as on $\bbz^d$ with $d\ge2$, $\rho_c<1/2$ (for the
continuous-time symmetric FEP on the ladder, $\rho_c\approx0.4755$ or,
for the slightly different dynamics of \cite{L}, $\rho_c\approx0.4874$
\cite{I}).  The results of \cite{I} also suggest that \cconj holds for
this model (although the scaling behavior of the $L_i(\delta)$ is not
discussed).  Rigorously, one may observe that for $\rho<\rho_c$ the
portions of the system to the left and right of an empty square are
independent under the $t\to\infty$ limiting measure $\nu_\rho$,
implying in particular that the locations of the empty squares (when
these have a nonzero density) form a renewal process and that the
portions of the system between them are jointly independent.  Further,
we see that the critical density must satisfy $\rho_c\ge1/4$, since at
smaller densities there would always be a finite density of empty
squares and the stationary state would be frozen \cite{I}. We have no
such lower bound for $\rho_c$ on $\bbz^d$, $d\ge2$.

 \medskip\noindent
  {\bf Acknowledgments:} We thank Cesar Ramirez-Ibanez for helpful
discussions, and two anonymous referees for very useful comments.

This version of the article has been accepted for publication, after
peer review, but is not the Version of Record and does not reflect
post-acceptance improvements, or any corrections. The Version of
Record is available online at:
http://dx.doi.org/10.1007/s00220-025-05441-z.

 \medskip\noindent
  {\bf Data availability statement:} This paper has no associated data.

 \medskip\noindent
  {\bf Funding and competing interests statement:} The authors have no
relevant financial or non-financial interests to disclose, and did not
receive support from any organization for the submitted work.

 \end{document}